\documentclass[final,5p,times,twocolumn]{elsarticle}

\usepackage{amssymb}
\usepackage{amsmath}
\usepackage{amsthm}
\usepackage{geometry}
\usepackage{multirow}
\newtheorem{theorem}{Theorem}[section]

\newtheorem{prop}{Proposition}[section]

\setcitestyle{square}

\makeatletter
\long\def\@makecaption#1#2{
  {\bfseries #1} #2
  \vskip\abovecaptionskip
  \vskip\belowcaptionskip}
\makeatother

\journal{Operations Research Letters}

\begin{document}

\begin{frontmatter}



\title{Bounding probability of small deviation on sum of independent random variables: Combination of moment approach and Berry-Esseen theorem}



\author{Jiayi Guo \fnref{1}}
\author{ Simai He \fnref{2}}
\author{ Zi Ling} 
\author{Yicheng Liu}

\fntext[1]{Supported by the Fundamental Research Funds for the Central Universities} 
\fntext[2]{Supported by NSFC Grants 71771141, 71825003}
\fntext[]{Email: guo.jiayi@sufe.edu.cn (J. Guo), simaihe@sufe.edu.cn (S. He), elaine\_lingzi@outlook.com (Z. Ling), sufelyc@163.com (Y. Liu)}

\address{School of Information Management and Engineering, Shanghai University of Finance and Economics, China}

\begin{abstract}


For small deviation bounds, i.e., the upper bound of probability  \text{Prob}$\left[\sum_{i=1}^n X_i \geq \sum_{i=1}^n \mathop{\mathbb{E}}[X_i]+\delta \right]$ where $\delta$ is small or even negative, many classical inequalities (say Markov's inequality, Chebyshev's inequality, Cantelli's inequality \cite{cantelli1910intorno}) yield only trivial, or non-sharp results, see (\ref{Feige}). In this particular context of small deviation, we introduce a common approach to substantially sharpen such inequality bounds by combining the semidefinite optimization approach of moments problem \cite{bertsimas2005optimal} and the Berry-Esseen theorem \cite{Shevtsova2010An}. As an application, we improve the lower bound of Feige's conjecture \cite{Feige04onsums} from 0.14 \cite{GARNETT2020105119} to 0.1798. 
\end{abstract}

\begin{keyword}
robust optimization \sep moment problem \sep sum of random variables  \sep probability of small deviation \sep  Berry-Esseen theorem 


\end{keyword}

\end{frontmatter}

\section{Introduction}
\label{}
The problem of upper bounding
\begin{equation} \label{sum version}
\text{Prob}\left[ \sum_{i=1}^n X_i \geq  \sum_{i=1}^n \mathop{\mathbb{E}}(X_i) + \delta \right].
\end{equation}
for independent random variables $X_i$ and a given constant $\delta>0$, has been studied for years. 
Many classic tail bounds of this type, such as  Markov's inequality, Chebyshev's inequality \cite{lin2011probability}, Hoeffding's inequality \cite{hoeffding1963}, Bennett's inequality \cite{bennett1962probability} and Bernstein's inequality \cite{bernstein1924modification} and applications have been well-studied in literature and textbooks \cite{boucheron2013concentration,lin2011probability}. However, those inequalities are designed when $\delta$ is large. For example, Hoeffding's inequality indicates the following:
\begin{equation*}
\text{Prob}\left[ \sum_{i=1}^n X_i \ge \sum_{i=1}^n \mathop{\mathbb{E}}(X_i) + \delta \right] \le e^{-\frac{2\delta^{2}}{\sum_{i=1}^n (b_i - a_i)^{2}}}.
\end{equation*}
where $X_i \in [a_i, b_i]$. If $\delta$ is a relatively small constant, this inequality only provides bounds that are not so sharp. In particular, when $\delta$ is $0$, it yields a trivial bound.

\begin{table*}
\caption{: Feige's bound in literature and our paper.}
\label{tab1} 
\centering
\setlength{\tabcolsep}{9mm}{
\begin{tabular}{cccc}
\hline
Name&Bound&Method&Theorem\\
\hline
Feige \cite{Feige04onsums}& $\frac{1}{13}$&\\
He \cite{he2010}&0.125& approximate MP(1,2,4) &\\
Garnett \cite{GARNETT2020105119}& 0.14 & approximate MP(1,2,3,4) & \\
\multirow{4}*{Our paper}&0.1536& approximate MP(1,2,4) and B-E& Theorem \ref{thm3}\\
&0.1541&  MP(1,2,4) and B-E & Theorem \ref{thm2}\\
&0.1587&  MP(1,2,3,4) and B-E& Theorem \ref{thm4}\\
&0.1798&  MP(1,2,3,4) and B-E with refinement& Theorem \ref{thm5}\\
\hline
\end{tabular}}
\end{table*}

In this context of small deviation, which is widely applied in graph theory \cite{Feige04onsums} and inventory management \cite{wang2015process}, there are limited general tools to derive such bound. One approach is to formulate this problem as a moments problem (MP). Given the moments information, we can further derive an equivalent semidefinite programing (SDP) problem to the original moments problem through duality theory \cite{bertsimas2005optimal} and sum-of-square technique \cite{Lasserre07},  based on the classical theorem established by the great mathematician Hilbert in 1888 as the univariate case of Hilbert's 17th problem, which states that an univariate polynomial is nonnegative if and only it can be represented as sum of squares of polynomials.

It is also worth to mention that Berry-Esseen theorem, a uniform bound between the cdf of sum of independent random variables and the cdf of standard normal distribution, is a quite powerful inequality, no matter $\delta$ is large or small. Specifically, we can define a random variable $S=\frac{\sum_{i=1}^n X_i}{\sqrt{\sum_{i=1}^n \mathop{\mathbb{E}}[X_i^2]}}$ with $\text{Var}[S]=1$ to represent the normalized sum, and a third moment bound $\psi_0=\frac{\sum_{i=1}^n \mathop{\mathbb{E}}[|X_i|^3]}{(\sum_{i=1}^n \mathop{\mathbb{E}}[X_i^2])^{3/2}}$. Then Berry-Esseen theorem indicates the following:
\begin{equation}
\sup_{x \in \mathbb{R}} |F(x)-\Phi(x)| \leq c_0 \psi_0
\end{equation}
where $F(x)$ and $\Phi(x)$ are the cdf of $S$ and standard normal distribution respectively with best known $c_0=0.56$ \cite{Shevtsova2010An}. Unsurprisingly, if all random variables are independent identically distributed, then central limit theorem states that their properly scaled sum tends to towards a normal distribution. Berry-Esseen theorem just provides a quantitative rate of convergence. 

Therefore, it is a nature idea to combine moment approach and Berry-Eseen theorem to achieve a better bound of small deviation problems.

As an application to better illustrate the way of combination, Feige\cite{Feige04onsums} first established a bound $\alpha=\frac{1}{13}$ and conjectured the true bound to be $\alpha=\frac{1}{e}$ of the following small deviation problem:
\begin{equation} \label{Feige}
\text{Prob}\left[ \sum_{i=1}^n X_i \geq 1 \right] \leq 1-\alpha > 0
\end{equation}
where $X_1,X_2,...,X_n$ are independent random variables, with $\mathop{\mathbb{E}}[X_i]=0$ and $X_i \geq -1$ for each $i$.  This inequality has many applications in the field of graph theory \cite{ferber2019uniformity}, combinatorics \cite{Alon2012}, and evolutionary algorithms \cite{Corus2014}. 

One contribution of this paper is to improve the Feige's bound from best-known $\alpha=0.14$ to $0.1798$ step by step as it is shown in  Table \ref{tab1}.



The other contribution of this paper is to introduce a general approach to bound probability of small deviation by merging the moment problem approach and the Berry-Esseen Theorem. Suppose we have a sequence of independent random variables $X_i$ and a constant $\delta$, and define $\sum_{i=1}^n \mathop{\mathbb{E}}(X_i^2)=D$. Here is the guideline of our approach.
\begin{itemize}
\item In order to achieve the upper bound of (\ref{sum version}), without loss of generality, we can assume $X_i$ has support set of at most $k$ discrete points where $k$ is the number of given moment information (including the trivial $0$-th order moment) on $\sum X_i$ by constructing an associated linear programming. This insight is an extension of lemma 6 in Feige \cite{Feige04onsums}, and has been established by Bertismas et. al. \cite{bertsimas2005optimal}.

\item For both the Berry-Esseen theorem and moment approach, it requires the distributions $X_i$ has bounded support,  i.e., there exists a certain constant $K$ such that $|X_i| \leq K$ for all $i$. When the distributions are not bounded, we divide the distributions into bounded and unbounded groups, and treats the unbounded group separately as in \cite{he2010}.

\item We can derive a bound of (\ref{sum version}) by Berry-Esseen theorem. Suppose we have $|X_i| \leq K$ for all $i$. Then $\psi_0 =\frac{\sum_{i=1}^n \mathop{\mathbb{E}}[|X_i^3|]}{D^{3/2}}\leq \frac{K}{D^{1/2}}$. When $D$ is large, it follows that $\psi_0$ is relatively small, and therefore Berry-Esseen theorem can provide a rather tight bound. 

\item We can also bound (\ref{sum version}) by the moment approach. In particular, when the distributions are bounded, we can bound the third moment or above through the second moment $D$. Theorem \ref{more} indicates the more moments we use, the better bound we can achieve. In addition, when $D$ is small, the bound of moment approach is often better as it is shown in theorem \ref{thmMono}.

\item We observe that often the worse-case scenario of these two approaches do not agree with each other, which provides us a great opportunity to merge these two methods together to further improve the bound estimation. In Section 3 and 4, we discuss how to synthetically merge these two approaches together to achieve better results. 

\end{itemize}

\section{An SDP formulation of the moment problem}
In this section, we introduce the classical SDP formulation of the moment problem.
Supposing $X$ is a real random variable and $P$ is a set of given moments, then we formulate moments problem as the following.
\begin{equation}\label{MP}
\begin{aligned}
    \max_X & &\text{Prob}[X \geq 0] & \\
     \text{subject to} & & \mathop{\mathbb{E}}[X^i]=M_i, & \text{ where } i \in P
\end{aligned}
\end{equation}
Note that $M_0=1$. The corresponding dual problem is
\begin{equation}\label{SPD}
\begin{aligned}
     \min_y \,\,\, & & \sum_{i \in P} y_i M_i & \\
     \text{subject to} & & \sum_{i \in P} y_i x^i \geq \text{\bf 1}_{x \geq 0}, & \,\,\,\text{ for all }  x \in R 
\end{aligned}
\end{equation}
where $\text{\bf 1}_{x \geq 0}$ is an indicator function.

This is an well-studied optimization problem and the dual formulation was first established in \cite{isii1960extrema} and treated extensively in \cite{bertsimas2005optimal}. In fact, the property of strong duality was shown in \cite[Theorem 2.2]{bertsimas2005optimal}. Moreover,
the dual constraint requires the polynomial function to be nonnegative,
which is equivalent to certain matrices being positive semidefinite, as it is shown in the following theorem. 

\begin{theorem}\cite[Section 3.a]{reznick2000some}\label{sos}

 \noindent A real polynomial function $f(x) =\sum_{r=0}^{2n} y_rx^r$ is nonnegative if and only if $f(x) = g(x)^2 + h(x)^2$ for some polynomial function $g(x)$ and $h(x)$. Furthermore, the nonnegativity of $f(x)$ implies the existence of an $(n + 1) \times (n + 1)$  positive semidefinite matrix $V$ such that $f(x) = X^TVX$ with $X = (1,x,x^2,...x^n)^{T}$. 

 \end{theorem}
 
 \noindent{\bf Proof.} 
 
 (SOS Decomposition): If $f(x)=g(x)^2+h(x)^2$, then it is obviously nonnegative. 
 
 If $f(x)$ is nonnegative, then all real roots of $f(x)$ are of even multipliers, because otherwise $f(x)$ will be negative locally. By the Fundamental Theorem of Algebra, $$f(x)=\prod_{j=1}^m (x - r_j)^{2n_j}\prod_{i=1}^{n-\sum_{j=1}^{m} n_j}(x - z_i)(x - \overline z_i).$$ Notice that $(x - z_i)(x - \overline z_i) = (x - a_i)^2 + b_i^2$, $f(x)$ can be decomposed as products of sum square of two functions. Since $(a^2 +b^2)(c^2+d^2) = (ac+bd)^2 + (ad-bc)^2$, then we can reformulate $f$ into $g^2 + h^2$.
 
 (SDP Representation): Suppose $f(x) =  g(x)^2 + h(x)^2$, and the coefficient vector of $g$ and $h$ are $u$ and $v$, respectively. Then $g(x) = (1,x,x^2,...x^n)u$, and $h(x) = (1,x,x^2,...x^n)v$. Let $V=uu^T + vv^T$, and we have $f(x) = g(x)^2 + h(x)^2 = X^Tuu^TX + X^Tvv^TX = X^T(uu^T+vv^T)X = X^TVX$.
 \qed
 
 Theorem 2.1 is the univariate case for Hilbert’s 17th problem. In fact, this theorem together with the further work by Lasserre \cite{Lasserre07} can help us transform the dual problem into a semidefinite program, as it is shown in the following proposition.
 
\begin{prop}\cite[Proposition 3.1]{bertsimas2005optimal} \label{Prop} 

\begin{itemize}
    \item The polynomial $g(x)=\sum_{r=0}^{2n}y_rx^r$ satisfies $g(x) \geq 0$ for all $x \in R$ \text{\bf if and only if} there exists a positive semidefinite matrix $V=[v]_{i,j=0,1,...,n}$ such that 
    $$y_r= \sum_{i,j:i+j=r}v_{ij}, \,\, r=0,...,2n $$
    \item The polynomial $g(x)=\sum_{r=0}^{n}y_rx^r$ satisfies $g(x) \geq 0$ for all $x \geq 0$ \text{\bf if and only if} there exists a positive semidefinite matrix $V=[v]_{i,j=0,1,...,n}$ such that 
    $$0= \sum_{i,j:i+j=2l-1}v_{ij}, \,\, l=1,...,n $$
    $$y_r= \sum_{i,j:i+j=2l}v_{ij}, \,\, l=0,...,n $$
\end{itemize}

\end{prop}

In all, moments problem (\ref{MP}) can be solved by its SDP formulation.

One key property of the moments problem is to achieve a better bound by taking advantage of additional moment information, as extra moment information yields a more restrictive constraint set in the moment problem.

\begin{theorem} \label{more}
 Given a real random variable $X$, consider the primal problem
 \begin{equation*}
 \begin{aligned}
   \text{\bf opt }(P)=\max_X & \,\,\,\,\, \text{Prob}[X \geq 0] & \\
     \text{subject to} & \,\,\,\,\, \mathop{\mathbb{E}}[X^i]=M_i, & \text{ where } i \in P
\end{aligned}
\end{equation*}
 Supposing we have $P_1 \subset P_2$, then $\text{\bf opt }(P_1) \geq \text{\bf opt }(P_2)$.
 \end{theorem}

Similarly, the following theorem explores the monotonicity of the moments problem with upper and lower bounds, as the feasible region enlarges as $D$ grows larger.


\begin{theorem} \label{thmMono}
Given a real random variable $X$ and mutually exclusive sets $P_1, P_2, P_3$, consider the primal problem where $B_i$ and $L_i$ are the upper and lower of moments in a function of a real number $D$. 
 \begin{align*}
   \text{\bf opt }(D)= \max_X & &\text{Prob}[X \geq 0] & \\
     \text{subject to} & & \mathop{\mathbb{E}}[X^i]=M_i, & \text{ where } i \in P_1 \\
     & & \mathop{\mathbb{E}}[X^i] \leq B_i(D), & \text{ where } i \in P_2 \\
     & & \mathop{\mathbb{E}}[X^i] \geq L_i(D), & \text{ where } i \in P_3
\end{align*}
 Supposing we have 
 \begin{itemize}
     \item $B_i(D) \geq 0$ and $B_i(D)$ is an increasing function in $D$ for all $i \in P_2$;
     \item $L_i(D) \leq 0$ and $L_i(D)$ is an decreasing function in $D$ for all $i \in P_3$,
 \end{itemize}
 then $\text{\bf opt }(D)$ is an monotonically increasing function in $D$.
\end{theorem}

\section{A combination of moment approach and Berry-Esseen theorem}
\begin{theorem} \label{MPBE}
 Let $X=\sum_{i=1}^n X_i$, and $D=Var(X)=\sum_{i=1}^n \mathop{\mathbb{E}}[ X_i]^2$ for independent random variables $X_i$ with bound $|X_i| \leq K$. In addition, suppose there exist mutually exclusive sets $P_1=\{0,1,2\}$, $P_2$, $P_3$ with increasing nonnegative functions $B_i(D)$ for $i \in P_2$ and decreasing nonpositive functions $L_i(D)$ for $i \in P_3$.
Then
$$\text{Prob}[X \geq 0] \leq \min_{D > 0}\max\{\text{\bf opt}(D), F(D)\}$$
where 
$$F(D)=0.5+0.56 \frac{K}{D^{1/2}}$$
and
 \begin{align*}
   \text{\bf opt }(D)= \max_X & \,\,\,\, \text{Prob}[X \geq 0] & \\
     \text{subject to} \,\,\,\, & \mathop{\mathbb{E}}[X^0]=1, & \\
    &   \mathop{\mathbb{E}}[X]=\sum_{i=1}^n \mathop{\mathbb{E}}[X_i], & \\
    &  \mathop{\mathbb{E}}[X^2]=D, & \\
     &  \mathop{\mathbb{E}}[X^i] \leq B_i(D), & \text{ where } i \in P_2 \\
     &  \mathop{\mathbb{E}}[X^i] \geq L_i(D), & \text{ where } i \in P_3
\end{align*}
Moreover, $\arg\min_{D > 0}\max\{\text{\bf opt}(D), F(D)\}$ is at the intersection of function $\text{\bf opt}(D)$ and function $F(D)$, if it exists.
\end{theorem}

\noindent{\bf Proof.} 
We can apply Berry-Esseen theorem on $\text{Prob}[X \geq 0]$.
\begin{align*}
    \text{Prob}[X \geq 0] & \leq 1 - \Phi(0) + c_0 \psi_0 \\
    & \leq 0.5 + 0.56 \frac{\sum_{i=1}^n |X_i|^3}{D^{3/2}} \\
    & \leq 0.5+0.56 \frac{K}{D^{1/2}} = F(D)
\end{align*}
When $D$ is large, Berry-Esseen theorem is effective, because $F(D)$ is a decreasing function. When $D$ is small, the moments problem performs well because $\text{\bf opt}(D)$ is an increasing function by theorem \ref{thmMono}.

Therefore, if $\text{\bf opt}(D)$ and $F(D)$ exists an intersection, then it is the optimal solution of $\min_{D > 0}\max\{\text{\bf opt}(D), F(D)\}$ due to the monotonicity of these two functions.
\qed

\section{Example: improve the bound of Feige's inequality}

In this section, we will show that a combination of moment approach and Berry-Essen theorem can improve Feige’s bound. 

As we see, Feige's conjecture (\ref{Feige}) has no assumptions on the upper bound of each random variable, though we know the lower bound is $-1$. The following theorem \ref{thm1} allows us to transform such variables $X_i$ into a group of corresponding $Y_i$ with both upper and lower bounds, through truncating the sufficiently large negative part of $X_i$ and rescale the rest. Similar technique was used in \cite{he2010, GARNETT2020105119}. In this way, we can apply the inequalities of sum of independent random variables with both upper bound and lower bound, as theorem \ref{MPBE} indicates.

\begin{theorem}\label{thm1}
Supposing for $m$ random variables $Y_1$, $Y_2$,...,$Y_m$  with mean zero and $ -\xi \leq Y_i \leq 1$ for some  fixed $0 < \xi \leq 1$, there exists an universal bound $\omega>0$ independent of $m$ and the choice of $Y_i$ such that 
$$\text{Prob}[\sum_{i=1}^m Y_i \leq \xi] \geq \omega$$. 

Then, consider $n$ random variables $X_1$, $X_2$,...,$X_n$  with mean zero and $X_i \geq -1$. 
$$\text{Prob}[\sum_{i=1}^n X_i <1] \geq e^{- {\xi} } \cdot \omega$$
\end{theorem}

\noindent{\bf Proof.}
 As it is shown in \cite{Feige04onsums}, without loss of generality, we can assume $X_i$ follows a two-point distribution.

Therefore, we can assume that there exists $0<a_i \leq 1$ and $b_i >0$ such that
$$X_i=\left\{
\begin{array}{c l}	
     -a_i  \text{ with probability } \frac{b_i}{a_i+b_i}\\
\\
     b_i  \text{ with probability } \frac{a_i}{a_i+b_i}
\end{array}\right.$$
given $E[X_i]=0$.

Suppose $b_1 \geq b_2 \geq ... \geq b_n$. Then we consider to make a partition and define $A=\{1,2,...,N\}$ and $B=\{N+1,...,n\}$ by a fixed number $\tau>1$ where
$$N=\max\{\,0,\,\max\{\,k \,|\, b_k \geq \tau (\sum_{i=1}^k a_i), 1\leq k \leq n\}\}$$
Define $a=\sum_{i=1}^k a_i$ and we have 
\[b_i \geq b_N \geq \tau a 
, \text{ for every } i \leq N\]
\[b_i \leq b_{N+1} \leq \tau(a+a_{N+1}) \leq \tau(a+1) 
, \text{ for every } i > N\]
If $N >0$, then
\begin{align*}
\text{Prob}[\sum_{i=1}^N X_i=-a ] &= \Pi_{i=1}^N \text{Prob}[X_i=-a_i]\\
&=\Pi_{i=1}^N (1-\frac{a_i}{a_i+b_i}) \geq  \Pi_{i=1}^N (1-\frac{a_i}{a_i+\tau a})\\ 
& \geq \Pi_{i=1}^N e^{-\frac{a_i}{\tau a}}=e^{-\frac{1}{\tau}}
\end{align*}
Therefore, 
\begin{align*}
\text{Prob}[\sum_{i=1}^n X_i < 1 ] & \geq \text{Prob}[\sum_{i \in A} X_i=-a ] \text{Prob}[\sum_{i \in B} X_i < a+1 ]\\
& \geq e^{-\frac{1}{\tau}} \text{Prob}[\sum_{i \in B} X_i < a+1 ]
\end{align*}

Let $Y_i=\frac{1}{\tau} \frac{X_i}{a+1}$ for $i \in B$. Note that $$Y_i =\frac{1}{\tau} \frac{X_i}{a+1} \leq   \frac{1}{\tau} \frac{\tau (a+1)}{a+1}=1$$ and 
$$Y_i =\frac{1}{\tau} \frac{X_i}{a+1} \geq -\frac{1}{\tau}.$$
Then, if we set $\xi=\frac{1}{\tau}$,
\begin{align*}
\text{Prob}[\sum_{i=1}^n X_i < 1 ] & \geq e^{-\frac{1}{\tau}} \text{Prob}[ \sum_{i \in B} X_i < a+1 ] \\
&= e^{-\frac{1}{\tau}} \text{Prob}[ \sum_{i \in B} Y_i < \frac{1}{\tau} ] \\
& \geq e^{- \xi} \cdot \omega
\end{align*}
\qed

For the rest of work, we will consider the following problem: Let $Y_1,Y_2,...,Y_n$ be independent random variable with mean zero and $-\xi \leq Y_i \leq 1$ for some $0< \xi \leq 1$. Without loss of any generality, we can assume
\begin{equation} \label{setupY}
Y_i=\left\{
\begin{array}{c l}	
     -a_i  \text{ with probability } \frac{b_i}{a_i+b_i}\\
\\
     b_i  \text{ with probability } \frac{a_i}{a_i+b_i}
\end{array}\right.
\end{equation}
where $0 \leq a_i \leq \xi$ and $0 \leq b_i \leq 1$. Then for any $n$, we are interested in the lower bound of 
$$\text{Prob}[\sum_{i=1}^n Y_i \leq \xi] , $$
as a key to improve Feige's bound $\alpha$ in (\ref{Feige}).

\subsection{Grouping the first, second and fourth moment information} \label{124}

When it comes to Feige's bound (\ref{Feige}), He and et al. improved it to $1/8$ by solving the moments problem with the first, second and fourth moment information. Therefore, we only consider the same moment information in this subsection as a fair comparison. 

Suppose we have $Y_1,Y_2,...,Y_n$ be independent random variables with mean zero and $-\xi \leq Y_i \leq 1$ for some $0< \xi \leq 1$, as it is stated in (\ref{setupY}). Let $Y=\sum_{i=1}^n Y_i$ and $D=Var(Y)=\sum_{i=1}^n E(Y_i^2)=\sum_{i=1}^n a_ib_i$.  Berry-Esseen theorem implies the following:
\begin{align*}
\text{Prob}[ Y \leq \xi ] &=\text{Prob}[ \frac{Y}{\sqrt{D}} \leq \frac{\xi}{\sqrt{D}} ] \\
& \geq \Phi(\frac{\xi}{\sqrt{D}}) - c_0 \psi_0 \\
& \geq \Phi(\frac{\xi}{\sqrt{D}}) - 0.56 \frac{\sum_{i=1}^n E[|Y_i^3|]}{D^{3/2}} \\
& \geq \Phi(\frac{\xi}{\sqrt{D}}) - 0.56 \frac{\max_i\{|Y_i|\} \sum_{i=1}^n E[|Y_i^2|]}{D^{3/2}} \\
& = \Phi(\frac{\xi}{\sqrt{D}}) - 0.56 \frac{1}{\sqrt{D}} 
\end{align*}
Define
\begin{equation} \label{F1}
    F_1(\xi,D)=\Phi(\frac{\xi}{\sqrt{D}}) - 0.56 \frac{1}{\sqrt{D}} 
\end{equation}
as a lower bound of $\text{Prob}[ Y \leq \xi ]$.

For moments problem, we can define $Z=Y-\xi$. Then:

\begin{itemize}
    \item $M_1=\mathop{\mathbb{E}}[Z]=-\xi$
    \item $M_2=\mathop{\mathbb{E}}[Z^2]=D+\xi^2$
    \item $M_4 = \mathop{\mathbb{E}}[Z^4] $
    \begin{align*}
        &=3D^2+6\xi^2D+\xi^4+\sum_{i=1}^n(\mathop{\mathbb{E}}[Y_i^4]-4\xi \mathop{\mathbb{E}}[Y_i^3]-3(\mathop{\mathbb{E}}[Y_i^2])^2)\\
        &=3D^2+6\xi^2D+\xi^4+\sum_{i=1}^n a_i b_i(a_i^2+b_i^2-4a_ib_i-4\xi(b_i-a_i))
    \end{align*}
\end{itemize}

Since $a_i^2+b_i^2-4a_ib_i-4\xi(b_i-a_i)$ is a convex function of $a_i$ when $b_i$ and $\xi$ are fixed, and is a convex function of $b_i$ when $a_i$ and $\xi$ are fixed. Then supposing we fix $\xi$, the optimal solution of 
$$\max_{0 \leq a_i \leq \xi, 0 \leq b_i \leq 1} a_i^2+b_i^2-4a_ib_i-4\xi(b_i-a_i)$$
is in the set  $\{(0,0), (\xi,0),(0,1),(\xi,1)\}$ with the optimal value $S(\xi)$. 
It follows that
$$M_4 \leq 3D^2+6\xi^2D+\xi^4+S(\xi) D$$

Let $\text{\bf opt}(\xi,D)$ be the optimal value of the following moments problem given $\xi$ and $D$.
\begin{equation*}
 \begin{aligned}
    \text{\bf opt}(\xi,D)=\max_X & \,\,\,\, \text{Prob}[Z \geq 0] & \\
     \text{subject to} \,\,\,\, & \mathop{\mathbb{E}}[Z^0]=1, & \\
    &   \mathop{\mathbb{E}}[Z]=-\xi, & \\
    &  \mathop{\mathbb{E}}[Z^2]=D+\xi^2, & \\
     &  \mathop{\mathbb{E}}[Z^4] \leq 3D^2+6\xi^2D+\xi^4+S(\xi) D, &
\end{aligned}
\end{equation*}

Define
\begin{equation} \label{F2}
F_2(\xi,D)=1- \text{\bf opt}(\xi,D)    
\end{equation}
 to be another lower bound of $\text{Prob}[ Y \leq \xi ]$. In addition, we can calculate $F_2(\xi,D)$ numerically by solving a corresponding SDP problem introduced in proposition \ref{Prop}, when $\xi$ and $D$ are fixed.

\begin{theorem}\label{thm2}
 Let $X_1$, $X_2$,...,$X_n$ be n independent random variables with $E[X_i]=0$ and $X_i \ge -1$ for each $i$ and let $X = \sum_{i=1}^n X_i$, then $$\text{Prob}[X \le 1] \ge 0.1541.$$
\end{theorem}

\noindent{\bf Proof.} 

Set $\xi=0.2$. Then $S(\xi) \leq 0.2$.
\begin{itemize}
    \item If $D \geq 2.374$, then 
    $$\text{Prob}[\sum_{i=1}^n X_i <1] \geq e^{-0.2}F_1(0.2, 2.374) >0.1541,$$
    as $F_1(0.2,D)$ is an increasing function in $D$.
    \item If $D \leq 2.374$, then  
    $$\text{Prob}[\sum_{i=1}^n X_i <1] \geq e^{-0.2}F_2(0.2, 2.374) > 0.1541,$$
    as $F_2(0,2,D)$ is an decreasing function in $D$ by theorem (\ref{thmMono}).
\end{itemize}

\begin{figure}[h!]
\centering
\caption{: Feige's bound by $F_1(0.2, D)$ and $F_2(0.2, D)$.}
\includegraphics[width=9cm]{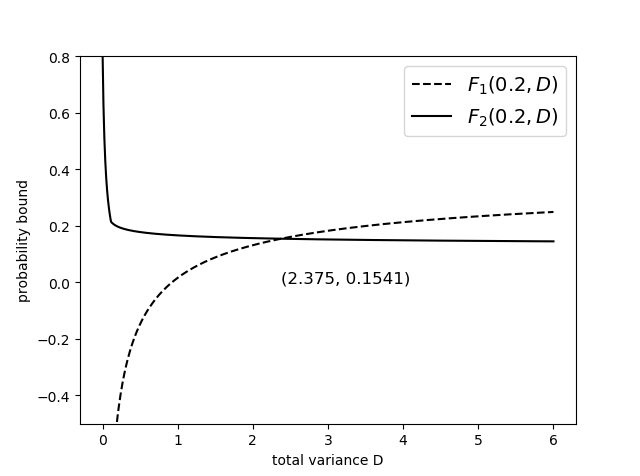}
\label{sosBE}
\end{figure}

In figure \ref{sosBE}, we plot $F_1(0.2, D)$ and $F_2(0.2,D)$ over the value of $D$. 

In all, the bound is improved to 0.1541. 
\qed

Instead of achieving bound 0.1541 numerically, we can roughly verify this result by an approximation of $F_2(\xi, D)$ in an explicit form. Specially, we can derive bound 0.1536 exactly as theorem \ref{thm3} indicates in the appendix.

\subsection{Add the third moment information} 
Recently, Garnett improved Feige's bound to 0.14 by a finer consideration of first four moments of the corresponding moments problem \cite{GARNETT2020105119}. If adding the third moment information, then we have the following lower bound in the same set-up as the previous section.
\begin{equation}\label{M3}
\begin{aligned}
     M_3 & = \mathop{\mathbb{E}}[Z^3]\\
     & = -\xi^3-3\xi D +\sum_{i=1}^n  \mathop{\mathbb{E}}[Y_i^3] \\
     & = -\xi^3-3\xi D - \sum_{i=1}^n a_i b_i ( a_i-b_i)\\
     & \geq -\xi^3-3\xi D -\xi D =  -\xi^3-4\xi D
\end{aligned}
\end{equation}

Let $\text{\bf opt}(\xi,D)$ be the optimal value of the following moments problem given $\xi$ and $D$.
\begin{equation*}
 \begin{aligned}
    \max_X & \,\,\,\, \text{\bf opt}(\xi,D)=\text{Prob}[Z \geq 0] & \\
     \text{subject to} \,\,\,\, & \mathop{\mathbb{E}}[Z^0]=1, & \\
    &   \mathop{\mathbb{E}}[Z]=-\xi, & \\
    &  \mathop{\mathbb{E}}[Z^2]=D+\xi^2, & \\
    &  \mathop{\mathbb{E}}[Z^3] \geq  -\xi^3-3\xi D -\xi D =  -\xi^3-4\xi D, & \\
     &  \mathop{\mathbb{E}}[Z^4] \leq 3D^2+6\xi^2D+\xi^4+S(\xi) D, & 
\end{aligned}
\end{equation*}

Define 
\begin{equation} \label{F4}
    F_4(\xi,D)=1- \text{\bf opt}(\xi,D)
\end{equation} 
to be another lower bound of $\text{Prob}[ Y \leq \xi ]$. 
Unsurprisingly, $F_4(\xi,D)$ should be better than $F_2(\xi,D)$.

\begin{theorem}\label{thm4}
 Let $X_1$, $X_2$,...,$X_n$ be n independent random variables with $E[X_i]=0$ and $X_i \ge -1$ for each $i$ and let $X = \sum_{i=1}^n X_i$. Then $$\text{Prob}[X \le 1] \ge 0.1587.$$
\end{theorem}

\noindent{\bf Proof.} 

Set $\xi=0.2$. Then $S(\xi) \leq 0.2$.
\begin{itemize}
    \item If $D \geq 2.464$, then 
    $$\text{Prob}[\sum_{i=1}^n X_i <1] \geq e^{-0.2}F_1(0.2, 2.464) >0.1587,$$
    as $F_1(0.2,D)$ is an increasing function in $D$.
    \item If $D \leq 2.464$, then  
    $$\text{Prob}[\sum_{i=1}^n X_i <1] \geq e^{-0.2}F_4(0.2, 2.464) > 0.1587,$$
    as $F_4(0.2,D)$ is an decreasing function in $D$ by theorem \ref{thmMono}.
\end{itemize}

\begin{figure}[h!]
\centering
\caption{: Feige's bound by $F_1(0.2, D)$ and $F_4(0.2,D)$.}
\includegraphics[width=9cm]{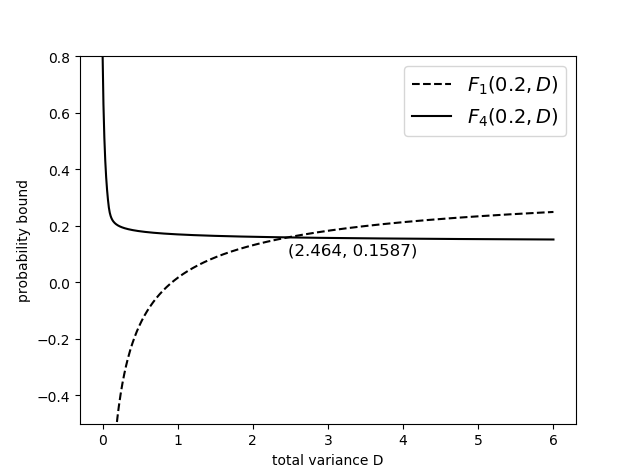}
\label{sosThird}
\end{figure}

In figure \ref{sosThird}, we plot $F_1(0.2, D)$ and $F_4(0.2,D)$ over the value of $D$. 
In all, the bound is improved to 0.1587. 
\qed

As we see, from 0.1541 to 0.1587, the Feige's bound was improved only a little. The reason is that we bound $M_3$ purely by $a_i - b_i \leq 1$ in (\ref{M3}), which is not enough. In fact, we observe that often the worse-case scenario of bound $M_3$ and the bound of Berry-Esseen term $\sum_{i=1}^n \mathop{\mathbb{E}}[|Y_i^3|]$ do not agree with each other. Therefore, we are able to better bound $M_3$ through the term $\sum_{i=1}^n \mathop{\mathbb{E}}[|Y_i^3|]$ (which is bounded as $\max{|Y_i|} \cdot D$ through this paper).
In general, this technique is significant to improve our result when we hybrid the moments method and Berry-Esseen theorem.

\begin{theorem}\label{thm5}
 Let $X_1$, $X_2$,...,$X_n$ be n independent random variables with $E[X_i]=0$ and $X_i \ge -1$ for each $i$ and let $X = \sum_{i=1}^n X_i$. Then $$\text{Prob}[X \le 1] \ge 0.1798.$$
\end{theorem}

\noindent{\bf Proof.} 

Define 
 $T_B=\sum_{i=1}^n \mathop{\mathbb{E}}[Y_i^3]=\sum_{i=1}^n a_i b_i \frac{a_i^2+b_i^2}{a_i+b_i}$
 in the same set-up as (\ref{setupY}). When applying Berry-Esseen theorem, we can define $\hat{F}_1(\xi,D,T_B)$ to be following
\begin{align*}
    \text{Prob}[\sum_{i=1}^n Y_i  \geq \xi] & \geq \Phi(\frac{\xi}{\sqrt{D}})-0.56\frac{\sum_{i=1}^n \mathop{\mathbb{E}}[Y_i^3] }{D^{3/2}} \\
    & = \Phi(\frac{\xi}{\sqrt{D}})-0.56\frac{T_B}{D^{3/2}} = \hat{F}_1(\xi,D,T_B).
\end{align*}

At the same time, define $T_M=\sum_{i=1}^n a_i b_i (a_i - b_i)$, and we have 
$$M_3=\mathop{\mathbb{E}}[Z^3]= -\xi^3-3\xi D - T_M.$$

Note that
$$T_B+T_M = \sum_{i=1}^n a_i b_i \frac{2a_i^2}{a_i+b_i} \leq \sum_{i=1}^n a_i b_i \cdot 2a_i \leq 2 \xi D $$

In this way, we can better bound the third moment $M_3$.

\begin{itemize}
    \item If $T_B=D$, then the Berry-Esseen bound remains the same i.e. $F_1(\xi,D)=\hat{F_1}(\xi,D,D)$. In this way, $T_M \leq (2\xi-1)D$ implying $M_3 \geq -\xi^3-3\xi D -(2\xi-1)D = -\xi^3-5\xi D +D$.
    \item If $\xi D \leq T_B \leq D$, then the Berry-Esseen bound improves i.e. $F_1(\xi,D) \leq \hat{F_1}(\xi,D,T_B)$. In this way, $T_M \leq 2\xi D- T_B$ implying $M_3 \geq -\xi^3-3\xi D -(2\xi D -T_B)=-\xi^3-5\xi D +T_B$.
    \item If $T_B < \xi D$, then the bound $T_B+T_M \leq 2 \xi D $ is no longer effective. We have $M_3 \geq -\xi^3-3\xi D -(2\xi D -T_B)=-\xi^3-4\xi D$, the same as the bound as (\ref{M3}) when $T_B=\xi D$.
\end{itemize}

For each given $T_B$, we can achieve corresponding bound of $M_3$ by the analysis above. Then we can define $\hat{F_4}(\xi,D,T_B)$ in a similar way.

Fix $\xi=0.2$. Suppose $T_B=s \cdot D$ for some $0 < s \leq 1$. Define function possible Feige's bound $g(s)$ to be the following: $$g(s) = \min_{D}\max\{e^{-0.2} \cdot \hat{F_1}(0.2,D,sD), \,\, e^{-0.2} \cdot \hat{F_4}(0.2,D,sD)\}$$ 

\begin{figure}[h!]
\centering
\caption{: Feige's bound under different s.}
\includegraphics[width=9cm]{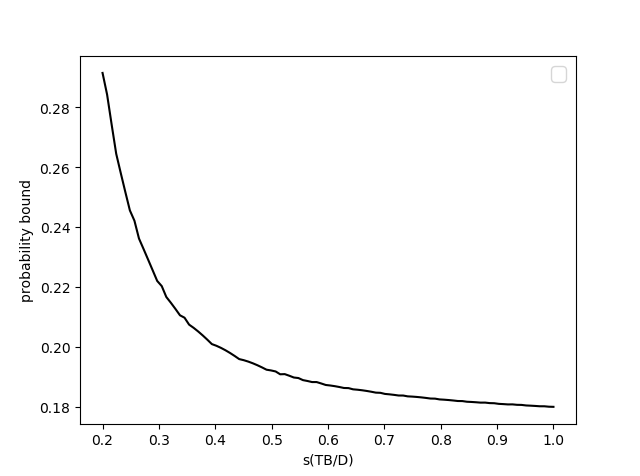}
\label{sosTB}
\end{figure}

Figure \ref{sosTB} plots the value of g(s) under different value of s.

Note that $\hat{F_1}(0.2,D,sD)$ is a decreasing function in $s$ for each given $D$, and $\hat{F_4}(0.2,D,sD)$ is an increasing function in $s$ for each given $D$. Figure \ref{sosTB} indicates the influence of improving Berry-Essen bound $\hat{F}_1$ dominates the influence of improving moment bound $\hat{F}_4$.

Therefore, we can set $T_B=D$.
\begin{itemize}
    \item If $D \geq 2.938$, then 
    $$\text{Prob}[\sum_{i=1}^n X_i <1] \geq e^{-0.2}\hat{F_1}(0.2, 2.938, 2.938) >0.1798,$$
    as $\hat{F}_1(0.2,D,D)$ is an increasing function in $D$.
    \item If $D \leq 2.938$, then  
    $$\text{Prob}[\sum_{i=1}^n X_i <1] \geq e^{-0.2}\hat{F}_4(0.2, 2.938,2.938) > 0.1798,$$
    as $\hat{F}_4(0.2, D,D)$ is an decreasing function in $D$.
\end{itemize}

In figure \ref{sosThird(Interplay))}, we plot $\hat F_1(0.2, D, D)$ and $\hat F_4(0.2,D, D)$ over the value of $D$. 

\begin{figure}[h!]
\centering
\caption{: Feige's bound by $\hat F_1(0.2, D, D)$ and $\hat F_4(0.2, D, D)$.}
\includegraphics[width=9cm]{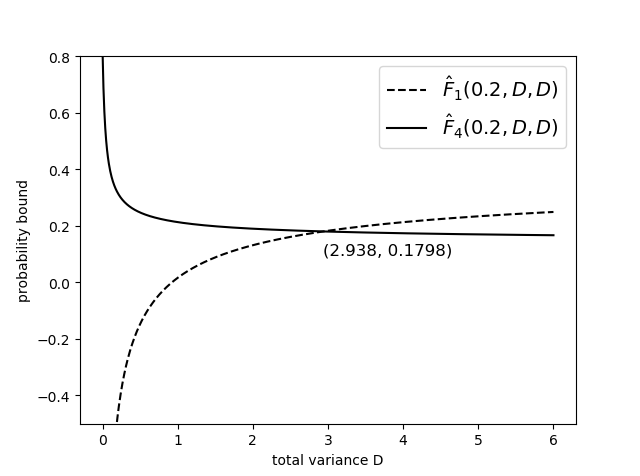}
\label{sosThird(Interplay))}
\end{figure}

In all, the bound is improved to 0.1798. 
\qed

\section{Summary}
\label{}
In this paper, we show that the combination of Berry-Esseen theorem and moment approach can better bound probability in small deviation. As an application, we improve Feige's bound from 0.14 to 0.1798 using first four moments. However, there is still a gap between 0.1798 to the conjectured $\frac{1}{e}$. Due to the length of this paper, we leave the readers to further improve it by including higher order moments, or better bounding fourth moment via $T_B$.

More importantly, we expect this common approach to be widely applied on other interesting small deviation problems. For example, Ben-Tal and et al. \cite{ben2002robust} conjectured the following: Consider a symmetric matrix $B \in R^{n\times n}$, and let $\xi=(\xi_1,\xi_2,...,\xi_n) \in R^n$ with coordinates $\xi_i$ of $\xi$ being independently identically distributed random variables with
$$Pr(\xi_i=1)=Pr(\xi_i=-1)=\frac{1}{2}.$$
Then,
$$Pr(\xi^T B\xi \leq Tr(B))\geq \frac{1}{4}.$$
Define the lower bound
$$y(n)=\inf_{B \in S^{n \times n}} Pr(\xi^T B\xi \leq Tr(B)).$$
The best known of result is $y(n) \geq \frac{3}{100}$ \cite{he2010}. Besides, Yuan showed the upper bound of $y(n)$ is $\frac{14}{64}$ by an example \cite{yuan2013counter}. Straightforward application of the  approach in this paper leads to improved bound of $\frac{6}{100}$ at least. Due to the space limitation, and since we believe finer consideration could vastly improve that bound, we omit the detailed proof and leave it for future research.


In all, it will be interesting to see our approach substantially sharpening inequality bound of small deviation problems and facilitating their applications. 

\appendix
\section{}
Supposing only the first, second and fourth moment information is considered as the set-up in section \ref{124}, we can solve the SDP exactly and achieve $F_2(\xi,D)$  to be the following:
$$\text{\bf opt}(\xi,D) =\max_{v} \frac{4(2\sqrt{3}-3)}{9} (-\frac{2M_1}{v}+\frac{3M_2}{v^2}-\frac{M_4}{v^4})$$ by \cite[Theorem 2.1]{he2010}.

Therefore, we can derive the following approximate bound to be $F_3(\xi,D)$ by choosing $v=\sqrt{2 M_4/3 M_2}$ for simplicity, as it was in \cite[Theorem 3.1]{he2010}. Thus,
\begin{small}
\begin{align*}
& \text{Prob}[\sum_{i=1}^n Y_i <\xi] \\
& \geq \frac{4(2\sqrt{3}-3)}{9} \inf_{D} {\bigg(}\sqrt{ \frac{6(D\xi^2+\xi^4)}{3D^2+(6+s)D\xi^2+\xi^4}}+\frac{\frac{9}{4}(D+\xi^2)^2}{3D^2+(6+s)D\xi^2+\xi^4} {\bigg)} \\
& = F_3(\xi, D)
\end{align*}
\end{small}
where $s=\max\{5, \frac{1}{\xi^2}-\frac{4}{\xi}, \frac{1}{\xi^2}-\frac{8}{\xi}+5 \}$.

As it is shown in figure \ref{sosFourth123}, when fixing $\xi=0.2$, the exact solution $F_2(0.2,D)$ is always above explicit approximate solution $F_3(0.2,D)$.

In addition, as a reproduction of Feige's bound derived in \cite{he2010}, $$\text{Prob}[X \leq 1] \geq e^{0.2} \cdot \lim_{D \rightarrow \infty} F_3(0.2,D) \geq 1/8$$
due to the monotonicity of $F_3(0.2,D)$.
\begin{figure}[h!]
\setcounter{figure}{0}
\caption{: Comparison of $F_2(0.2,D)$ and $F_3(0.2,D)$.}
\includegraphics[width=9cm]{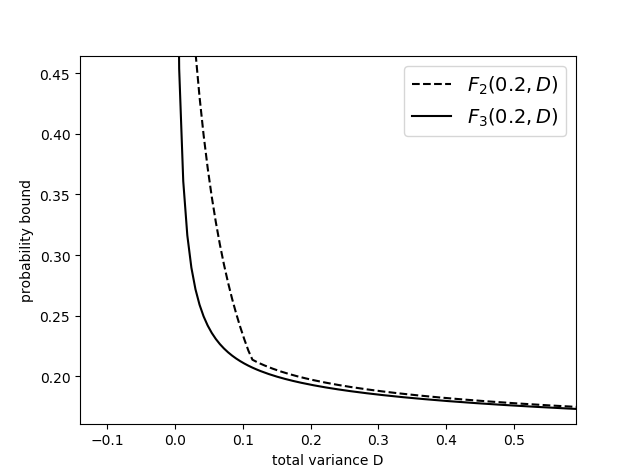}
\label{sosFourth123}
\end{figure}

If combining $F_1$ and $F_3$, we will have the following theorem.

\begin{theorem}\label{thm3}
 Let $X_1$, $X_2$,...,$X_n$ be n independent random variables with $E[X_i]=0$ and $X_i \ge -1$ for each $i$ and let $X = \sum_{i=1}^n X_i$, then $$\text{Prob}[X \le 1] \ge 0.1536.$$
\end{theorem}

\noindent{\bf Proof.} 

If we set $\xi=0.2$ as we discussed above.
\begin{itemize}
    \item If $D \geq 2.367$, then 
    $$\text{Prob}[\sum_{i=1}^n X_i <1] \geq e^{-0.2}F_1(0.2, 2.367) >0.1536,$$
    as $F_1(0.2,D)$ is an increasing function in $D$.
    \item If $D \leq 2.367$, then  
    $$\text{Prob}[\sum_{i=1}^n X_i <1] \geq e^{-0.2}F_3(0.2, 2.367) > 0.1536,$$
    as $F_3(0.2,D)$ is an decreasing function in $D$.
\end{itemize}



In all, Feige's bound is 0.1536. 
\qed

  \bibliographystyle{elsarticle-num} 
  \bibliography{Bounding_Probability_of_small_deviation}






\end{document}